\documentclass[12 pt]{amsart}
\usepackage{tikz}
\theoremstyle{plain}
\newtheorem{theorem}{Theorem}
\newtheorem{corollary}{Corollary}
\newtheorem{lemma}{Lemma}
\newtheorem{proposition}{Proposition}
\usepackage{color}
\theoremstyle{definition}

\theoremstyle{remark}

\newtheorem*{remark}{Remark}

\newcommand{\A}{\alpha}
\newcommand{\Z}{\mathbb{Z}}
\newcommand{\N}{\mathbb{N}}
\newcommand{\R}{\mathbb{R}}
\newcommand{\Q}{\mathbb{Q}}
\newcommand{\E}{\epsilon}

\title{On the Frequency of Balanced Times in Cylinder Flows}

\author{Jon Chaika}
\address{Jon Chaika\\ Rice University Department of Mathematics\\ 6100 S. Main St. \\ Houston, TX 77005}
\email{jmc5125@rice.edu}

\author{David Ralston}
\address{David Ralston \\ OSU Department of Mathematics \\ 231 West 18th Ave. \\ Columbus, OH 43210}
\email{ralston@math.ohio-state.edu}


\date{\today}

\begin{document}
\begin{abstract}
Given an irrational $\A \in (0,1)$ and an $x \in [0,1)$, the set of \emph{balanced times}, for which the same number of $k\A+x \mod{1}$ are less than or equal to one half as are larger than one half, is in general infinite, but sparse in terms of density.  We investigate the sparseness of this sequence in terms of summation over reciprocals.  Our results are that for the generic pair $(\A,x)$, the resulting sum diverges,  but there are certain exceptional $\A$ for which the associated sums converge for every $x$.
\end{abstract}

\maketitle

\section{Introduction and Survey of Results}

Consider the simple irrational rotation of the unit circle: $X = \R / \Z = [0,1)$, coupled with the Lebesgue measure-preserving transformation $R_{\A}(x) = x+ \A \mod{1}$ (we generally omit the ``mod 1" from our discussion, and $\A \notin \Q$ is always assumed).  If we define $I=[0,.5]$, results dating back at least to Weyl show that we asymptotically expect roughly as many points in the sequence $\{x+n\A\}$ to be in $I$ as not: if we define \[f(x) = \chi_{I}(x) - \chi_{x \setminus I}(x),\] then
\[\lim_{n \rightarrow \infty} \frac{1}{n} \sum_{i=0}^{n-1}f(x+i\A)=0.\]  

However, the ergodic sums \[\sum_{i=0}^{n-1} f(x+i\A)\] are unbounded for all $x$ \cite{kesten}, so these sums have arbitrarily large deviations from their expected value; we are interested in how often the ergodic sums are \emph{exactly} zero.  Following \cite{schmidt}, define the \emph{cylinder flow} $T_{\A}: X \times \Z \rightarrow X \times \Z$ by \[T_{\A}(x,n)=\left(x+\A,n+f(x)\right).\]  The projections onto the first and second coordinates are given by $\pi_1(x,n)=x$ and $\pi_2(x,n)=n$, respectively referred to as the \emph{position} and \emph{level} of the pair $(x,n)$.

In order to investigate how often the ergodic sums for $x$ achieve their limiting average in the case of the circle rotation, we are interested in how often the pair $(x,0)$ returns to the zero-level in the cylinder flow.  So, we define the set of \emph{balanced times} for $x$ by \[V_{\A,x}=\left\{n \in \Z^+ : \pi_2 \left( T_{\A}^n (x,0) \right)=0\right\},\] and $V_{\A,x}(N)=V_{\A,x}\cap\left\{1,2,\ldots,N\right\}$.  

The transformation $T_{\A}$ is ergodic if and only if $\A$ is irrational \cite{aaronson-keane} \cite{conze-keane} \cite{oren}, so for almost every $x$, the set $V_{\A,x}$ is infinite, but of density zero: \[\lim_{N \rightarrow \infty} \frac{1}{N} \#(V_{\A,x}(N)) = 0,\]  where $\#(S)$ denotes the cardinality of a set $S$.  The purpose of this paper is to consider the size of the set $V_{\A,x}$ in terms of the divergence or convergence of the sums \[\sum_{n \in V_{\A,x}} \frac{1}{n}.\]  The transformation $T_{\A}$ allows us to consider times when the ergodic sums of $f(x)$ exactly zero as a recurrence property of the pair $(x,0)$ in the
cylinder flow.

We will consistently analyze the transformation $T_{\A}$ through the continued fraction expansion of $\A$.  We refer the reader to \cite{khinchin} for an excellent treatment of continued fractions, pausing here only to give the briefest of summaries:

For each irrational $\A \in (0,1)$, there is a unique sequence of $a_i \in \N$ so that \[\A= \lim_{n \rightarrow \infty} \frac{p_n}{q_n},\] where \[\frac{p_n}{q_n}=\cfrac{1}{a_1+\cfrac{1}{a_2+\cfrac{1}{\ddots+\frac{1}{a_n}}}}.\]  The fractions $p_n/q_n$ are called \emph{convergents} for $\A$, and the $a_i$ are called the \emph{partial quotients} of $\A$.  Both $p_n$ and $q_n$ obey the same recursion formula in terms of the partial quotients:
\[p_{n+1}=a_{n+1}p_n+p_{n-1},\] \[q_{n+1}=a_{n+1}q_n + q_{n-1},\] where $p_0/q_0=0/1$ and $p_1/q_1=1/a_1$.  We write $\A=[a_1,a_2,\ldots]$.

The function $\gamma:[0,1) \rightarrow [0,1)$ given by truncation of the continued fraction expansion of $\A$ by removing the first entry is called the \emph{Gauss map}: \[\gamma([a_1,a_2,\ldots,a_n])=[a_2,a_3,\ldots,a_n] \textrm{  for }\A \in \Q,\] \[\gamma([a_1,a_2,\ldots])=[a_2,a_3,\ldots] \textrm{ for }\A \notin \Q,\] where by convention $\gamma([a_1])=\gamma(0)=0$.  This transformation is ergodic with respect to the \emph{Gauss measure} on $X$ (which we denote by $G$), a measure which is mutually absolutely continuous with respect to Lebesgue measure (we denote Lebesgue measure by $L$).  For any Lebesgue-measurable set $S$, \[\frac{1}{2\log 2}L(S) \leq G(S) \leq \frac{1}{\log 2} L(S).\]  When considering some variable $\A$, we write $a_i(\A)$ for the $i$-th partial quotient of $\A$.  If $\A$ is fixed, however, we generally write simply $a_i$.  We also define the quantities \[A_n(\A) = \sum_{i=1}^n a_i(\A),\] where again, if $\A$ is fixed, we write only $A_n$.

Our central result is that the general situation is one where the set $V_{\A,x}$ is fairly `large':

\begin{theorem}\label{theorem-almost sure divergence}
{For Lebesgue almost every pair $(\A,x)$, \[\sum_{n \in V_{\A,x}}\frac{1}{n}=\infty.\]}
\end{theorem}

However, Theorem \ref{theorem-almost sure divergence} cannot be strengthened to a statement about \emph{every} $x$: for all $\A$ there is a unique point $x_0$ for which $V_{\A,x_0}=\emptyset$ \cite{ralston - rotations}.  More remarkable, however, is the fact that Theorem \ref{theorem-almost sure divergence} can not be strengthened to a statement about every $\alpha$.

\begin{theorem}\label{theorem-heaviness for convergence}
{There is an uncountable collection of $\A$ so that for \emph{every} $x$, \[\sum_{n \in V_{\A,x}} \frac{1}{n}<\infty.\] }
\end{theorem}

\section{Preliminary Results in Continued Fractions}\label{section - cont frac stuff}

For a sequence $S \subset \N$, the \emph{upper density} of $S$ is given by \[d^*(S) = \limsup_{n \rightarrow \infty} \frac{\#\left(S \cap \left\{1,2,3,\ldots,n\right\}\right)}{n}.\]

\begin{proposition}\label{proposition - summation}
{Suppose $b_n > b_{n+1}>0$, so that \[\sum_{n=1}^{\infty}b_n = \infty.\]  Then \[\sum_{n \in S}b_n = \infty\] for every $S \subset \N$ with $d^*(S)>0$ if and only if \[\liminf_{n \rightarrow \infty}n b_n >0.\]}
\begin{proof}

Assume that $(\liminf nb_n) > \delta>0$, and let $0<\E<d^*(S)$ for some $S$.  Let $n_i$ be an increasing sequence of integers so that \[\frac{\# \left( S \cap \{n_i+1,n_i +2, \ldots, n_{i+1} \} \right)}{n_{i+1}-n_i} > \E,\] as well as $n_{i+1}>2n_i$ (and $n_0=0$).  Then we see that 
\begin{align*}
\sum_{i \in S \cap \{1,2,\ldots,n_k\} } b_{i} &> \E\left(n_1 b_{n_1} + \left( n_2 -n_1\right) b_{n_2} + \ldots + \left(n_k - n_{k-1}\right) b_{n_k}\right) \\
&> \E \left( \frac{1}{2} n_1 b_{n_1} + \frac{1}{2} n_2 b_{n_2} + \ldots + \frac{1}{2} n_k b_{n_k}\right)\\
&> \frac{\E}{2} k \delta
\end{align*}
and therefore \[\sum_{i \in S} b_i = \infty.\]

Conversely, assume that $n_k$ are chosen so that $n_k b_{n_k} < 2^{-k}$, and $n_{k+1}>2^k n_k$.  Fix some $\E>0$, and construct the set 
\[S =\bigcup_{i=1}^{\infty} \{n_i, n_i +1, \ldots ,n_i + \left[\E n_i\right]\} \] where $[x]$ is the integer part of $x$.  Then $d^*(S) \geq (\E/2) >0$, and

\[ \sum_{i \in S} b_i = \sum_{i=1}^{\infty} \sum_{j=0}^{[\E n_i]} b_{n_i+j} 
< \sum_{i=1}^{\infty} \E (n_i- n_{i-1}) b_{n_i} 
< \sum_{i=1}^{\infty} \E n_i b_{n_i}
< \E. \]
\end{proof}
\end{proposition}

\begin{lemma}\label{lemma - pos measure upper density}
{Let $\{n_t\}$ be an increasing sequence of natural numbers ($t=1,2,\ldots$), and for any $\A$ define \[N_{\A} = \left\{ m : A_{n_m}(\A) < 12n_m \log(n_m) \right\}.\]  Then the set \[B=\left\{ \A : d^*(N_{\A}) > 0\right\}\] is of positive measure.}
\begin{proof}
Define \[S_r = \left\{ \A: a_i(\A)< n_r^2 \text{ for all } i<n_r\right\}.\]   We refer to a standard estimate in continued fractions \cite[\S 12]{khinchin}: for any choice of $a_1,a_2,\ldots,a_{n-1}$, if we define $B$ to be the set of those $\A$ whose continued fraction expansion begins with $a_1,a_2,\ldots,a_{n-1}$, then 
\begin{align*}
\frac{G \left( B \cap \left\{\A: a_n(\A)=k\right\}\right)}{G\left( B\right)} &\leq 2 G \left(\left\{\A: a_n(\A)=k\right\}\right) \\
&= 2G\left( \left\{\A: a_1(\A) = k\right\}\right)\\
&\leq \frac{2}{\log 2} L\left( \left\{\A: a_1(\A)=k\right\}\right)\\
&= \frac{2}{\log 2} L \left(\left[\frac{1}{k},\frac{1}{k+1}\right)\right)\\
&< \frac{2}{k^2 \log 2}
\end{align*}

\begin{align*}
\int_{S_r} A_{n_r}(\A) dG(\A) &< \sum_{j=1}^{n_r}\left( \sum_{i=1}^{n_r^2}i \frac{2}{i^2 \log 2}\right)\\
&< \frac{2}{\log 2} n_r \left(\log (n_r^2)+1\right)\\
&< 6 n_r \left( \log n_r +1\right)
\end{align*}
(the index $j$ tracks which partial quotient we are considering, and $i$ ranges over the allowed values for $a_j(\A)$).  Define \[C_r = \left\{ \A  \in S_r: A_{n_r} > 12 n_r (\log n_r+1)\right\}.\]  By Chebychev's inequality, then, $G \left(C_r \right) < 1/2$, from which it follows that we cannot have almost every point avoid $X \setminus C_r$ for a density-one collection of $r$, as the sets $S_r$ satisfy \[\lim_{r \rightarrow \infty} G(S_r) = 1,\] by \cite[Theorem 32]{khinchin}.
\end{proof}
\end{lemma}
\begin{remark} If $n_i=i$ then for almost every $\A$, there are infinitely many $r$ for which $\A \notin
C_r$.
\end{remark}
We may now prove:

\begin{theorem}
\label{theorem - cont frac theorem}
{For almost every $\A \in [0,1)$, \[\sum_{i=1}^{\infty} \frac{1}{A_i(\A)} = \infty.\]}
\begin{proof}
Let $n_t=2^t$ and $b_i(\A)=1/(A_{2^t}(\A))$ for all $2^{t-1}<i \leq 2^t$.  Note that $b_i(\A)<1/(A_i(\A))$, and \[\sum_{i=1}^{\infty}\frac{2^{t-1}}{A_{2^t}(\A)} = \sum_{i=1}^{\infty}b_i(\A).\]  By Lemma \ref{lemma - pos measure upper density}, we have a positive measure set of $\A$ for which there is a sequence of $t$, of positive upper density, along which \[\frac{2^{t-1}}{A_{2^t}(\A)} \geq \frac{2^{t-1}}{12(2^t \log(2^t))} = \frac{1}{24t \log(2)}.\]  We therefore see that \[\liminf_{t \rightarrow \infty} t b_t(\A)>0\] on this positive measure set of $\A$, so by applying Proposition \ref{proposition - summation}, the set
\[B= \left\{ \A: \sum_{i=1}^{\infty} \frac{1}{A_i(\A)}=\infty \right\}\] is of positive Lebesgue measure, so it is therefore of positive Gauss measure.  The set $B$ is clearly invariant under the Gauss map $\gamma$, which is ergodic with respect to the Gauss measure.  Therefore, $B$ is of full Gauss measure, and therefore of full Lebesgue measure.
\end{proof}
\end{theorem}

\section{Proof of Theorem \ref{theorem-almost sure divergence}}

It suffices to prove that for almost all $\A$, the set of $x$ for which the sum of reciprocals of elements of $V_{\A,x}$ diverges is of positive measure.  By ergodicity of $T_{\A}$ for $\A \notin \Q$, the induced transformation $T'_{\A}$ on $X \times\{0\}$, the transformation which sends a position $x$ to the position of its first return to the zero level under $T_{\A}$, is also ergodic.  It is clear that the sum of reciprocals of $V_{\A,x}$ diverges if and only if the sum of reciprocals of $V_{\A,T'_{\A} x}$ diverges.  Thus, the set of positions corresponding to divergent sums for a particular $A \notin \Q$ is invariant under an ergodic map.

\begin{lemma}
\label{lemma - max level by qn}
{The quantities $A_n$ dominate the maximum level reached by an orbit through time $q_n$: $\forall x \in[0,1), \A \notin \Q$, 
\[\max_{i=1,\ldots,q_n} \left|\pi_2 (T_{\A}^i(x,0)) \right| \leq 3 A_n(\A).\]
}
\begin{proof}
The distribution of the positions $x+ i \A$ for $i=1,2,\ldots,q_n$ is very regular, with one point in each interval of the form \[ [x+j/q_n,x+(j+1)/q_n],\] for $j=0,1,\ldots,q_n-1$ \cite[Theorem 1]{kesten}.  By considering where the positions $0$ and $1/2$ may fall relative to these intervals whether $q_n$ is even or odd, as well as whether the intervals containing $0$ and $1/2$ each have a point in the orbit which hits or misses $I$, \[\left|\pi_2\left( T_{\A}^{q_n}(x,0)\right)\right| \leq 3\] for all $n$.  We also remark that
\begin{equation}\label{equation - split sums}\pi_2 \left( T_{\A}^n(x,0)\right) = \pi_2 \left(T_{\A}^j(x,0) \right) + \pi_2 \left(T^{n-j}_{\A}(x+j\A,0) \right).\end{equation} 

We now prove the lemma inductively, noting that it is trivial for $n=0$ (in which case $A_0(\A)=0$ by default) or $n=1$.  Recall that $q_n=a_nq_{n-1}+q_{n-2}$.  The maximum level achieved by $T^i_{\A}(x,0)$ for $i=1,2,\ldots,q_{n-1}$ is bounded by $3A_{n-1}$ by our inductive hypothesis, and now at time $q_{n-1}$, the level is no larger than three in absolute value.  So by segmenting off portions of length $q_{n-1}$, we have:

\begin{center} \begin{tabular}{|c| c|}
\hline  Range of $i$ & Bound on level \\
\hline \hline $1,2,\ldots q_{n-1}$ & $3A_{n-1}$\\
\hline $q_{n-1}+1,\ldots,2q_{n-1} $& $3A_{n-1}+3$\\
\hline $2q_{n-1}+1,\ldots, 3q_{n-1}$ & $3A_{n-1}+6$\\
\hline \vdots & \vdots \\
\hline $(a_n-1)q_{n-1}+1,\ldots, a_nq_{n-1}$ & $3A_{n-1}+3(a_n-1)$\\
\hline $a_nq_{n-1}+1,\ldots,a_nq_n+q_{n-2}=q_n$ & $3A_{n-2}+3a_n$ \\
\hline
\end{tabular}
\end{center}
The lemma now follows from the fact that \[\max\left\{3A_{n-2}+3a_n, 3A_{n-1}+3(a_n-1)\right\} \leq 3A_n. \qedhere\]
\end{proof}
\end{lemma}

\begin{lemma}\label{messylemma}
{There is a positive constant $C>0$ so that for every $n$ (and fixed $\A$), \[L \left\{x \in X : \frac{1}{q_n}\#\left(V_{\A,x}(q_n)\right) > \frac{1}{4(6A_n +1)} \right\} \geq C.\]}
\begin{proof}
Partition $X$ into intervals whose endpoints are given by \[\{0, -\A, -2\A, \ldots, -(q_n-1)\A\}\] as well as by \[\left\{ \frac{1}{2}, \frac{1}{2}-\A, \frac{1}{2}-2\A, \ldots, \frac{1}{2}-(q_n-1)\A\right\}.\]  Then for $x,y$ in the same interval of this partition, for any $0\leq m \leq q_n$, we have \[\pi_2 \left(T_{\A}^m(x,0) \right) = \pi_2 \left(T_{\A}^m(x,0) \right).\]  Similarly, if $J$ is a partition element and $J+k\A$ is the interval obtained by translating $J$ by $k\A$ ($k\in\N$), then points in $J+k\A$ have identical levels for $m=1,2,\ldots,q_n-k$.

Now, for any $x \in X$, define \[t_x(N) = \left\{1 \leq m \leq N : \pi_2 \left( T_{\A}^m(x,0)\right) =t\right\}.\]  By Eq\eqref{equation - split sums} we see that if $j,k \in t_x(N)$ and $j>k$, then $j-k \in V_{\A,x+k\A}$.  So, for $k \in t_x(q_n)$:
\[ \# \left(V_{\A,x+k\A}(q_n) \right) \geq \# \left(t_x(q_n) \cap \left\{k+1,k+2,\ldots,q_n\right\} \right).\]  By Lemma \ref{lemma - max level by qn}, we have \[\bigcup_{t=-3A_n}^{3A_n} t_x(q_n) = \{1,2,\ldots,q_n\}.\]  If we define \[B= \left\{ t: \#\left(t_x(q_n)\right) > \frac{q_n}{2(6A_n+1)}\right\},\] then we have \[\# \left(\bigcup_{B} t_x(q_n) \right) \geq \frac{q_n}{2}\] by applying a pigeonhole argument.

Fix $t \in B$, and note that the least third of those $k \in t_x(q_n)$ have the property that \[\#\left(t_x(q_n) \cap \left\{k+1,k+2,\ldots,q_n\right\}\right) \geq \frac{1}{2}\#\left(t_x(q_n)\right) > \frac{q_n}{4(6A_n+1)}\] (of course, the quantity of such $k$ is exactly one half the total if $t_x(q_n)$ has an even number of elements, and one fewer if it is odd; we use one highly inaccurate but succinct and sufficient estimate for both cases).

So at least half of the points in the length $q_n$ orbit of an arbitrary point $x$ belong to one of $6A_n+1$ possible levels which is hit frequently, and at least a third of the returns to that level have the property that they have many future returns to their given level within $q_n$ steps.  Clearly, then, a positive proportion $C$ of the positions $x$ have the property that they have many returns to the zero level within the first $q_n$ steps, as quantified in the statement of the lemma (we need not spend much effort estimating $C$ except to note that it is positive and independent of $n$).  
\end{proof}
\end{lemma}

Let \[f(x,n)=\min\left\{\underset{i\in V_{\alpha,x}}{\overset{q_n}{\sum }} \frac 1 i, \underset{i}{\overset{n}{\sum }} \frac 1 {q_i} \frac {q_i} {4(6A_n+1)}\right\}.\]
It follows from Lemma \ref{messylemma} that for each $n$  a positive measure set of $x$ have $f(x,n)$ proportional to $\sum_{i=1}^{n} \frac 1 {A_i}$.  So by Theorem \ref{theorem - cont frac theorem}, for almost every $\A$ there is a set of $x$ of positive measure on which \[\limsup_{n \rightarrow \infty} f(x,n) =\infty.\]  That is, for almost every $\A$, there is a positive measure set of $x$ for which our sums diverge, and by ergodicity, for such $\A$, almost every $x$ gives rise to divergent sums.

\begin{remark}
In fact, we have a bit more information that simple divergence of our sums.  For almost every $\A$, the sums of $(A_n(\A))^{-1}$ diverge at least as fast as the sums of $(n \log n)^{-1}$ along a subsequence, and for such $\A$, almost every $x$ has the sums diverging at the same rate (along a possibly different subsequence). 
\end{remark}

\section{Proof of Theorem \ref{theorem-heaviness for convergence}}

The function $f(x)=\chi_{I}(x)-\chi_{X \setminus I}(x)$ is an upper-semicontinuous function on $X$ (recall that $I$ is closed), and the transformation $R_{\A}(x)=x+\A \mod{1}$ preserves Lebesgue measure, so the \emph{heavy set} \[H= \left\{ x \in X : \sum_{i=0}^{n-1} f(R_{\A}^i (x)) \geq 0, \hspace{.1 in} n=1,2,3,\ldots\right\}\] is nonempty \cite{peres}.  The relation to our cylinder flow is apparent: by our definition of $T_{\A}(x,n)$, we have \[\pi_2\left( T^n(x,0)\right) = \sum_{i=0}^{n-1}f( R_{\A}^i(x)).\]  In fact, in \cite{ralston - rotations} it is shown that for all irrational $\A$ there is a unique point $x_0$ so that \[\pi_{2}\left(T_{\A}^n(x_0,0) \right) =\sum_{i=0}^{n-1}f(R_{\A}^i(x_0)) \geq 1\] for every $n>1$.  Clearly, for this $x_0$ we have $V_{\A,x_0}=\emptyset$, and trivially \[\sum_{i \in V_{\A,x_0}} \frac{1}{i} < \infty.\]  However, using the techniques developed in \cite{ralston - rotations}, we may determine the set $V_{\A,\A}$ (note that we fix our initial position to be equal to our rotation amount) for a very specific class of $\A$, for which $V_{\A,\A}$ will be infinite but highly controllable: \[\A=[2a_1,b_1,2a_2,b_2,2a_3,b_3, \ldots].\]  We summarize the technique for constructing the orbit of $\A$ under rotation by $\A$ in Figure \ref{figure - heavy rotation pictures}.  The process is a direct result of explicit consideration of the first return to the interval $[0,1-2a_1\A)$, which induces a rotation by $[2a_2,b_2,2a_3,b_3,\ldots]$, indicating a clear inductive procedure for constructing the orbit.

\begin{figure}[b t]
\parbox{2.3 in}{
\centering{\begin{tikzpicture}[xscale=1.5][yscale=2]
\draw (0,0) -- (2,0);
\draw[->] (0,0) -- (.9,1);
\draw[->] (1,1) -- (2,.1);
\draw (1,2) node[anchor=south]{A single peak};
\draw (1,1.5) node[anchor=south]{of height $a_1$.};
\end{tikzpicture}}
}
\hfill
\parbox{2.3 in}{
\centering{\begin{tikzpicture}[yscale=1.5, xscale=.5]
\label{figure - first drawing}
\draw (0,0) -- (9,0);
\draw (0,0) -- (1,1) --  (2,0) -> (3,1) -> (4,0);
\draw[dashed][->] (4,0) -- (5,1);
\draw[dashed] (6,1)--(7,0);
\draw[->] (7,0) -- (8,1) -- (9,.1);
\draw (4.5,1.5) node[anchor=south]{Each peak height $a_1$,};
\draw (4.5,1.25) node[anchor=south]{$b_1$ total peaks.};
\end{tikzpicture}}
} 
\caption{To produce the levels through time $q_1=2a_1$, simply climb to level $q_1$ and then fall.  This single peak is then repeated $b_1$ times to obtain the picture through time $q_2-1$.}
\vspace{.5 in}
\parbox{2.3 in}{
\begin{tikzpicture}[xscale=.6, yscale=.8]
\draw (0,0) --  (9,0);
\draw (0,0) rectangle (1,1);
\draw (1,.4) rectangle (2,1.4);
\draw (3,3) rectangle (4,4);
\draw (4,3.4) rectangle (5,4.4);
\draw (5,3) rectangle (6,4);
\draw (7,.4) rectangle (8,1.4);
\draw (8,0) rectangle (9,1);
\draw [->] (.1,1.1) .. controls (.33,1.3) and (.66, 1.6)  .. (1,1.5) node[anchor=south] {+1};
\draw [->] (.1,1.1) .. controls (1,3) and (2,5)  .. (4,4.5) node[anchor=south] {+$a_{n+1}$};
\draw (5,2.1) node[anchor=south]{A single peak of height};
\draw (5,1.5) node[anchor=south]{  $a_1+a_2+\ldots+a_{n+1}$.};
\end{tikzpicture}
}
\hfill
\parbox{2.3 in}{
\centering{\begin{tikzpicture}[yscale=2, xscale=.5]
\label{figure - first drawing}
\draw (0,0) -- (9,0);
\draw (0,0) rectangle (2,.7);
\draw (2,0) rectangle (4,.7);
\draw[dashed][->] (0,0) -- (.2,.1) --(.4,0) -- (.8,.3) -- (1,.2) -- (1.2,.5) -- (1.4,.4) -- (1.6,.6) -- (1.8,.4) -- (2,0) -- (2.2,.1) --(2.4,0) -- (2.8,.3) -- (3,.2) -- (3.2,.5) -- (3.4,.4) -- (3.6,.6) -- (3.8,.4) -- (4,0) ;
\draw[dashed][->] (4.1,.5) -- (6.9,.5);
\draw[dashed][->] (7,0) -- (7.2,.1) --(7.4,0) -- (7.8,.3) -- (8,.2) -- (8.2,.5) -- (8.4,.4) -- (8.6,.6) -- (8.8,.4) -- (9,0);
\draw (7,0) rectangle (9,.7);
\draw (4.5,1.25) node[anchor=south]{Each peak same height,}; \draw (4.5,1) node[anchor=south]{$b_{n+1}$ total peaks.};
\end{tikzpicture}}
}
\caption{The general process is much the same; copy the figure through time $q_{2n}$ and move up one step at a time (the first figure contains an extra small peak, for a box of length $q_{2n}+q_{2n-1}$) a total of $a_{n+1}$ times, then falling back down (for a total length of $2a_{n+1}q_n+q_{n-1}=q_{n+1}$).  This new peak is copied $b_{n+1}$ times, once with an extra $q_{2n}$ peak, so that the total length is $q_{2n+2}$.}
\label{figure - heavy rotation pictures}
\end{figure}

Given the orbit through time $q_{2n}-1$, we construct the orbit through time $q_{2n+1}$ by translating the levels reached up by one a total of $a_n$ times, then back down.  One translate is extended by an addition length of $q_{2n-1}$, so that the total length achieved is $q_{2n+1}=2a_{n}q_{2n} + q_{2n-1}$.  This `single peak' is then copied $b_n$ times, once with an additional piece of size $q_{2n}$, to create the orbit of length $q_{2n+2}=b_nq_{2n+1}+q_{2n}$.  We refer to this process as ``stacking and shifting."  Our goal is to utilize large values of $a_i$ and small values of $b_i$ to tightly control how often any level is hit in the orbit of $\A$, and then transfer this knowledge to all $x$ via the unique ergodicity of the simple irrational circle rotation $R_{\A}$.

\begin{lemma}
Through time $q_{2n}$, exactly $A_n+1$ different levels are hit by the orbit of $(\A,0)$.
\begin{proof}
The proof is inductive, and clearly true for $n=0$ ($q_0=1$ and $A_0=0$).  We also know that the orbit of $\A$ through time $q_{2n}$ consists of a single peak, and our inductive hypothesis tells us that this peak hits all levels from $0$ to $A_{n}$.  To construct the orbit through time $q_{2n+2}$, this peak is copied $2a_{n+1}$ times, translated up $a_{n+1}$ times and then down.  So, this new single peak achieves a height of $A_n+a_{n+1}=A_{n+1}$, for a total of $A_{n+1}+1$ total levels hit.
\end{proof}
\end{lemma}

\begin{corollary}
For \emph{any} $x$, a $q_{2n}$-length portion of the orbit of $x$ hits at least $A_n$ different levels.
\begin{proof}
Regardless of which $x$ we consider, the $q_{2n}$-length orbit will be identical to the $q_{2n}$-length orbit of some $k\A$.  So, we need only show that any $q_{2n}$-length portion of the orbit of $\A$ hits $A_n$ different levels.  The initial segment of length $q_{2n}$ is a single peak of height $A_n$, which is then going to be shifted and stacked to form longer orbits.  So, regardless of where we begin looking in the orbit of $\A$, a portion of length $q_{2n}$ must contain either a complete `rise' or `fall' corresponding to one of these shifted/stacked peaks of length $q_{2n}$, which has a total height of $A_n$.
\end{proof}
\end{corollary}

When building the orbit of $\A$ through time $q_{2n}$, the peak formed through length $q_{2n-2}$ is stacked up $a_n$ times, then down, and then this new picture is shifted $b_n$ times.  Provided that the new peaks are relatively large compared to the old ones, then, if we look at the distribution of the number of hits to each of the $A_n$ levels which are achieved, the distribution will be fairly uniform:

\begin{lemma}\label{good lemma}
If $a_n>\E A_{n-1}$ for some fixed $\E>0$, then there is some $C=C(\A)<\infty$ so that for all $n$, the initial peak of length $q_{2n}$ (which is of height $A_n$) contains no more than $C(q_{2n}/A_n)$ hits to any level.
\end{lemma}

As before, knowledge about one orbit tells us knowledge about any orbit:

\begin{corollary} \label{q2n count}
For \emph{any} $x$, the number of hits to the zero level through time $q_{2n}$ is no larger than $2C(q_{2n}/A_n)$.
\begin{proof}
The $q_{2n}$-length orbit of any $x$ will occur somewhere in the orbit of $\A$.  As such, it will overlap at most two of the $q_{2n}$-length peaks, each of which has no more than $C(q_{2n}/A_n)$ hits to any given level.
\end{proof}
\end{corollary}

We may now complete the proof of Theorem \ref{theorem-heaviness for convergence}.  Choose any particular sequence $c_i$ so that \[\sum_{i=1}^{\infty} c_i < \infty.\]  We now fix  $\A=[2a_1,b_1,2a_2,b_2,\ldots]$.  Choose the partial quotients so that the $b_i$ are bounded by $M$ and the $a_i$ satisfy the conditions that there exists $r$ such that $(a_{i+1}+1)<A_i^r$ and 
\[ \frac{\log(2MA_n)}{A_n}<c_n .\]

Clearly, there are uncountably many such $\A$, and once $\A$ has been chosen, the constant $C$ from Lemma \ref{good lemma} is fixed. Consider \[\sum_{n \in V_{x,\alpha}}\frac{1}{n}=\sum_{k=0}^{\infty} \sum_{q_{2k}<n \in V_{x,\A} }^{q_{2k+2}} \frac {1} {n}\]

Breaking up this summation further,

\[\sum_{q_{2k}<n \in V_{x,\A}}^{q_{2k+2}} \frac {1} {n} \leq \sum_{q_{2k}<n \in V_{x,\A}}^{(a_{k+1}+1)(b_{k+1}+1)q_{2k}} \frac {1} {n}  = \sum_{i=1}^{(a_{k+1}+1)(2M)} \sum_{ iq_{2k}<n \in V_{x,\A}} ^{(i+1)q_{2k}} \frac {1} {n}\]

By Corollary \ref{q2n count} this is no larger than \begin{align*}
\sum_{i=1}^{(a_{k+1}+1)(2K)} \frac{2C(q_{2k})}{iq_{2k}A_k} &\leq \frac{2C}{A_k} (\log (a_{k+1}+1)2M)+1) \\
&\leq \frac{2Cr}{A_k} (\log (A_{k})2M)+1).\end{align*}

It follows that $\sum_{n \in V_{x,\alpha}}\frac{1}{n}$ is uniformly bounded for $x \in X$.

The proof can be refined (by observing that at most $2MA_k+3$ of the blocks $[iq_{2k},(i+1)q_{2k}]$ intersect $V_{x,\A}$) so that the condition $a_{i+1}<A_i^r$ is unnecessary.

\begin{remark}  In considering the proof of Theorem \ref{theorem-heaviness for convergence}, one should not lose the forest for the trees, so to speak.  The specific growth condition on the $a_i$ is not so important.  Rather, by forcing the $a_i$ to grow exceptionally quickly, we will be able for any $f(n)$ for which both $f(n) \rightarrow 0$ and $nf(n) \rightarrow \infty$ to construct some $\A$ so that
\[\lim_{n \rightarrow \infty}f(n) \max_{k \leq n} \left(\pi_2 (T^k(\A,0)) \right) = \infty.\]  With the maximal levels growing like $A_n$ (whose growth we may control through our choice of $a_i$),  the distribution of how many times each level has been hit is very unusual, allowing us to force the remarkable behavior that we seek.  The set of such $\A$ is clearly uncountable, being unchanged by slight (bounded, for example) perturbation of the $a_i$ and $b_i$.  Compare to \cite[Corollary 1]{kakutani-petersen}.  For an investigation of how this technique may be used to very precisely control the growth of the maximal height function, see \cite{ralston-arxiv}.
\end{remark}

The statement of Theorem \ref{theorem-heaviness for convergence} may actually be strengthened somewhat.  Let $\delta$ be some number strictly between $1/2$ and $1$, and consider the sum \[\sum_{n \in V_{x,\A}}\frac{1}{n^{\delta}}.\]  Construct $\A=[2a_1,b_1,\ldots]$, where the $b_i$ are again bounded by $M$, and having chosen some cenvergent series $c_n$, choose the $a_i$ so that 
\[\frac{q_{2n}^{1-\delta}}{A_n^{\delta}}<c_n.\]

We may be certain of building such an $\A$ (again, uncountably many) because $q_{2n} < (2a_n+1)(b_n+1)q_{2n-2}<3a_n(M+1)q_{2n-2}$.  We then have \[\frac{q_{2n}^{1-\delta}}{\left(a_n+A_{n-1}\right)^{\delta}} < \frac{\left(3(M+1)q_{2n-2} \right)^{1-\delta}}{a_n^{2\delta-1}}.\]  As $q_{2n-2}$ does not depend on our choice of $a_n$ and $\delta>1/2$, we may select $a_n$ so large so as to achieve our desired bound.

In this case, the relevant computation is given by
\begin{align*}
\sum_{n \in V_n}\frac{1}{n^{\delta}} &< \frac{2Cq_{2n}}{A_n q_{2n}^{\delta}} + \sum_{i=1}^{2MA_n+1}\frac{2Cq_{2n}}{A_n(iq_{2n})^{\delta}}\\
&=\frac{2Cq_{2n}^{1-\delta}}{A_n}\left( 1+ \sum_{i=1}^{2MA_n+1}\frac{1}{i^{\delta}}\right)\\
&< \frac{2Cq_{2n}^{1-\delta}}{A_n}\left( 1+\frac{(2MA_n)^{1-\delta}}{1-\delta}\right)\\
&< \frac{4C(M q_{2n})^{1-\delta}}{(1-\delta)A_n^{\delta}}\\
&< \frac{4CM^{1-\delta}}{1-\delta}c_n,
\end{align*}

where we have used the fact that \[\sum_{i=1}^N \frac{1}{i^{\delta}} < 1+ \int_1^{N-1} \frac{1}{x^{\delta}}dx.\]

So we can in fact claim that for any $\delta>1/2$, there is an uncountable collection of $\A$ for which every $x$ has the property that \[\sum_{n \in V_{x,\A}} \frac{1}{n^{\delta}}<\infty.\]

\section*{Acknowledgments}
The authors wish to thank William Veech for posing the problem studied herein and suggesting valuable references, as well as Michael Boshernitzan and Barak Weiss for stimulating conversations and encouragement.

\end{document}